\documentclass[11pt]{article}
\usepackage[dvips]{graphicx}
\usepackage{color,amsfonts,amssymb}
\usepackage{amsfonts,epsf,amsmath,tikz}
\usepackage{latexsym}
\usepackage{color}
\usepackage[ruled,vlined,linesnumbered]{algorithm2e}

\newtheorem{theorem}{Theorem}[section]
\newtheorem{lemma}[theorem]{Lemma}
\newtheorem{proposition}[theorem]{Proposition}

\newtheorem{property}[theorem]{Property}

\newcommand{\proof}{\noindent{\bf Proof.\ }}
\newcommand{\qed}{\hfill $\square$ \bigskip}
\newcommand{\gt}{{\rm gt}}
\newcommand{\SM}{{\rm SM}}

\begin{document}
\title{The geodesic-transversal problem}

\author{
	Paul Manuel$^{a}$
	\and
	Bo\v{s}tjan Bre\v{s}ar$^{b,c}$
	\and
	Sandi Klav\v zar$^{b,c,d}$
}

\date{}

\maketitle
\vspace{-0.8 cm}
\begin{center}
	$^a$ Department of Information Science, College of Computing Science and Engineering, Kuwait University, Kuwait \\
	{\tt pauldmanuel@gmail.com}\\
	\medskip
	
	$^b$ Faculty of Natural Sciences and Mathematics, University of Maribor, Slovenia\\
	{\tt bostjan.bresar@um.si}\\
	\medskip

	$^c$ Institute of Mathematics, Physics and Mechanics, Ljubljana, Slovenia\\
	\medskip

	$^d$ Faculty of Mathematics and Physics, University of Ljubljana, Slovenia\\
	{\tt sandi.klavzar@fmf.uni-lj.si}\\
\end{center}

\begin{abstract}
A maximal geodesic in a graph is a geodesic (alias shortest path) which is not a subpath of a longer geodesic. The geodesic-transversal problem in a graph $G$ is introduced as the task to find a smallest set $S$ of vertices of $G$ such that each maximal geodesic has at least one vertex in $S$. The minimum cardinality of such a set is the geodesic-transversal number $\gt(G)$ of $G$. It is proved that $\gt(G) = 1$ if and only if $G$ is a subdivided star and that the geodesic-transversal problem is NP-complete. Fast algorithms to determine the geodesic-transversal number of trees and of spread cactus graphs are designed, respectively.    
 \end{abstract}


\medskip
\noindent{\bf AMS Subj.\ Class.}: 05C69;  05C85; 68R10

\section{Introduction}
Given a set $U$ and a family $S = \{S_1, \ldots, S_k\}$, where $S_i\subseteq U$, a subset $H$ of $U$ is a {\em hitting set} for the family $S$ if $H \cap S_i$ $\neq$ $\emptyset$ for all $i\in\{1, \ldots, k\}$. The hitting set problem is to find a smallest hitting set for $S$. The hitting set problem is NP-complete \cite{GaJo79} and has been studied in different terminologies. In particular, in graph theory the term {\em $S$-transversal problem} presents the quest for a minimum set of vertices that intersect every set of a given family $S$ of subsets of the vertex set. When $S$ is a collection of maximal cliques of a graph, the $S$-transversal problem is called the {\em clique-transversal problem}~\cite{AnSc91, BaNa96, ChCh96, CoFo88, DaKr95, DaKr97, GuRa00}, and when $S$ is a collection of fixed size cliques, it is called the {\em generalized clique transversal problem}~\cite{CoFo88, DaKr95}. The clique-transversal problem is polynomially solvable for interval graphs and NP-complete for chordal graphs~\cite{DaKr95}.  Dahlhaus et al.~\cite{DaKr97} have studied the $S$-transversal problem where $S$ is a collection of hyperedges in a hypergraph. When $S$ is a collection of $k$-paths, the $S$-transversal problem is called the {\em $k$-path-transversal problem}. This problem has been well-studied under different terminologies~\cite{bjk13, bkk11, FuNu16, HoSu19, Manu18, MoMa20, TuYa13, Yann81}. 

A {\em geodesic} in a graph $G$ is a shortest path between two vertices, and  a geodesic is {\em maximal} if it is not a subpath of a longer geodesic. When $S$ is a collection of {\bf maximal geodesics}, we call the $S$-transversal problem the {\em geodesic-transversal problem}. A geodesic on $k$ vertices is a {\em $k$-geodesic}. When $S$ is a collection of $k$-geodesics, the $S$-transversal problem is called {\em $k$-geodesic-transversal problem}. 

To our knowledge, there is no literature on the geodesic-transversal problem and the $k$-geodesic-transversal problem. The {\em geodesic-transversal number} of $G$, denoted by $\gt(G)$, is the minimum cardinality of a geodesic-transversal set of $G$. A set $S$ of vertices is a {\em gt-set} of $G$ if $S$ is a minimum cardinality geodesic-transversal set of $G$. Thus, the geodesic-transversal problem of $G$ is to find a gt-set of $G$. It is easy to see that the $2$-geodesic-transversal problem is the vertex cover problem.

In the next section, we provide further motivation for the new geodesic-transversal problem. In Section~\ref{sec:basic}, we determine the geodesic-transversal number of some graphs and show that this number equals $1$ precisely for subdivided stars. We also prove that the geodesic-transversal problem is NP-complete for general graphs. In Section~\ref{sec:trees} we derive a  polynomial algorithm for arbitrary trees, while in Section~\ref{sec:cactus} a fast algorithm is designed for spread cactus graphs.  

%%%%%%%%%%%%%%%%%%%%%%%%%%%
%% M O T I V A T I O N %%%%
%%%%%%%%%%%%%%%%%%%%%%%%%%%

\section{Motivation from (large-scale) network theory}
\label{sec:motivation}

The geodesic-transversal problem is not entirely new. The path version of this problem is quite popular in graph theory and is well studied by graph theory researchers \cite{bjk13, bkk11, FuNu16, HoSu19, Manu18, MoMa20, TuYa13, Yann81}. A set $S$ of vertices of a graph $G$ is a {\em $k$-path vertex cover} if every path of order $k$ in $G$ contains at least one vertex from $S$~\cite{bkk11}. It is not uncommon in graph theory that the same concept is studied under different names. If indeed so, this indicates that the concept is of wider interest. The $k$-path vertex cover has been studied also as vertex $k$-path cover~\cite{bjk13}, $k$-path vertex cover~\cite{brkb-2017, bkk11, HoSu19, Manu18, MoMa20}, $VCP_k$-set~\cite{TuYa13}, and $k$-path cover~\cite{FuNu16}. The {\em $k$-path vertex cover problem} is to find the minimum cardinality of a $k$-path vertex cover. The problem is NP-hard for cubic planar graphs of girth 3~\cite{bkk11, TuYa13} and for bipartite graphs~\cite{Yann81}. The problem has applications in many areas, such as traffic control~\cite{TuZh11} and wireless sensor networks~\cite{bkk11}. Funke et al.~\cite{FuNu16} have provided a list of applications of this problem on different domains. The concepts of path transversal have also been generalized to the context of hypergraphs~\cite{zak}. The geodesic-transversal problem is a natural extension and adaptation of the path-transversal problem. Note that the $k$-path vertex cover problem and the $k$-geodesic transversal problem coincide in general graphs when $k=2$, and coincide in triangle-free graphs when $k=3$.

Betweenness centrality and closeness centrality are key measures of large-scale network analysis \cite{PeGe16, XiHa16}. The concepts of betweenness centrality and closeness centrality play a vital role in the study of large-scale network analysis including social networks \cite{GoFi13, LiLi15, PeGe16}, brain networks \cite{FoZa16, JoLa10, KwCh19}, biological networks (gene regulatory networks, protein-protein interaction network) \cite{KoSc08, KwCh19}, chemical networks \cite{ZhNa15}, communication networks \cite{ChSi13}, transport networks \cite{LiLi19, Ren15} and IoT networks \cite{Ren15, XiHa16} etc. 
The betweeness centrality $B(v)$ and closeness centrality $C(v)$ are defined as follows \cite{PeGe16, XiHa16}: 
\[ B(v) = \sum_{s\neq v \neq t} \frac{\sigma_{st}(v)}{\sigma_{st}} \]
\[ C(v) = \sum_{s\neq v \neq t} \sigma_{st}(v) \]
where $\sigma_{st}$ is the total number of geodesics from node $s$ to node $t$ and $\sigma_{st}(v)$ is the number of those paths that pass-through $v$. 

The scope of geodesic-transversal is wider than betweenness centrality and closeness centrality. The geodesic load geo-load($v$) of a vertex $v$ of a graph $G$ is defined as the number of maximal geodesics which traverse through $v$. The concept of geo-load of a network is applied in the geodesics-based routing algorithms \cite{Rack02, StKa17}. The concept is also used in load-balanced routing of fixed interconnection networks \cite{Sing05, XuHu19}. While the betweenness centrality of a vertex focuses on all possible geodesics, the geodesic load of a vertex concentrates on only maximal geodesics.

Some interesting combinatorial problems of large-scale network analysis are propagation (malware propagation \cite{YuGu15}, immunization \cite{Prak12}, disease propagation \cite{WaZh16} and data communication \cite{HoBr95}), broadcasting, and gossiping problems \cite{Xu02}. An interesting research problem is to demonstrate how the geodesic-transversal is a good model to represent these problems in large-scale network analysis. 

%%%%%%%%%%%%%%%%%%%%%%%%%%%
%%%%%  B A S I C S %%%%%%%%
%%%%%%%%%%%%%%%%%%%%%%%%%%%
\section{Basic observations and NP-completeness}
\label{sec:basic}

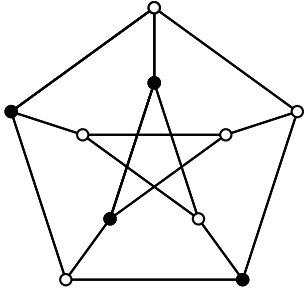
\begin{figure}[ht!]
\begin{center}
 \begin{tikzpicture}[scale=0.5,style=thick]
 \foreach \f in {0,72,...,360}
    \foreach \r in {2,4}
       { \draw (\f+90:4) -- (\f+90+72:4);
       \draw (\f+90:4) -- (\f+90:2);
       \draw (\f+90:2) -- (\f+90+72*2:2);}
\foreach \f in {0,72,...,360}
    \foreach \r in {2,4}
       {\draw[fill=white] (\f+90:\r) circle (0.15cm);}
 \foreach \f in {0,144}
    \foreach \r in {2}
   \draw[fill=black] (\f+90:\r) circle (0.15cm);
 \foreach \f in {72,216}
    \foreach \r in {4}
   \draw[fill=black] (\f+90:\r) circle (0.15cm);
  \end{tikzpicture}
  \caption{A gt-set of the Petersen graph}
\label{fig:Petersen}
\end{center}
\end{figure}

For a starting example consider the Petersen graph $P$. It is of diameter $2$, therefore to hit all the five maximal geodesics on the outer $5$-cycle we need at least two vertices. Similarly, we need at least two vertices to hit the maximal geodesics which are subpaths of the inner $5$-cycle. Hence $\gt(P)  \ge 4$. On the other hand, in Fig.~\ref{fig:Petersen} a geodesic-transversal set with four vertices is shown, hence we conclude that $\gt(P) = 4$. Using a similar reasoning we can deduce that if $r,s\ge 1$, then $\gt(K_{r,s})=\min\{r,s\}$. 

The following simple lemma will turn out to be quite useful. 

\begin{lemma}
\label{lem:equal-distance}
Let $Q$ be a geodesic of a graph $G$ and $x,y\in V(Q)$. If $u$ is a vertex from $V(G)\setminus V(Q)$ such that $d(u,x) = d(u,y)$, then $Q$ does not extend to a geodesic that contains $u$. 
\end{lemma}

\proof
Suppose on the contrary that $Q'$ is a geodesic such that $Q$ is contained in $Q'$ and $u\in V(Q')$. Clearly, on the geodesic $Q'$, the vertex $u$ cannot lie between $x$ and $y$. Therefore, either $d(u,x) < d(u,y)$ or $d(u,y) < d(u,x)$, and each of the possibilities in a contradiction with the lemma assumption. 
\qed

Clearly, $\gt(P_n)=1$ for all $n\in \mathbb{N}$. In particular, $\gt(P_n)=1$ because its only vertex forms a geodesic by itself and hence has to lie in its unique gt-set. Considering an arbitrary edge $e$ of the complete graph $K_n$, $n\ge 3$,  and a vertex not on the edge, Lemma~\ref{lem:equal-distance} implies that at least one of the endpoint of $e$ must lie in a geodesic-transversal set of $K_n$. Consequently, $\gt(K_n) = n-1$ holds for $n\ge 2$. These two examples generalize as follows, where by a {\em subdivided star} we mean the graph obtained from $K_{1,k}$, $k\ge 1$, by subdividing each of the edges of $K_{1,k}$ arbitrary number of times (possibly zero). If $k=1$, then the subdivided stars coincide with  the family of paths. 

\begin{proposition}
If $G$ is a connected graph of order at least $2$, then $1\le \gt(G)\le n(G)-1$. In addition, the lower bound is attained if and only if $G$ is a subdivided star, and the upper bound is attained if and only if $G$ is a complete graph of order at least $2$. 
\end{proposition}

\proof
Since every graph $G$ has at least one maximal geodesic, we infer $\gt(G)\ge 1$. Since every maximal geodesic of a non-trivial graph contains at least two vertices, we infer $\gt(G)\le n(G)-1$. 

Suppose now that $\gt(G) = 1$ and let $\{u\}$ be a gt-set of $G$. Let $T$ be a BFS-tree of $G$ with the root $u$. 

We first claim that $G$ is bipartite. Suppose on the contrary that there exists an edge $xy$ of $G$, where vertices $x$ and $y$ lie in the $k^{\rm th}$ distance level of $T$, for some $k\ge 1$. Then $d_G(u,x) = d_G(u,y) = k$.  Consider now an arbitrary maximal geodesic $Q$ of $G$ that contains the edge $xy$. Then Lemma~\ref{lem:equal-distance} implies, that $u$ does not belong to $Q$, a contradiction with the assumption that $u$ forms a gt-set. Hence the claim. 

We next claim that $G$ is a tree. Suppose on the contrary that $G$ contains at least one cycle $C$.  Since we already know that $G$ is bipartite, considering the cycle $C$ we infer that there exist a vertex $x$ of $C$ which lies in some $k^{\rm th}$ distance level of $T$ such that $x$ has two neighbors (in $G$), say $y$ and $z$, in the $(k-1)^{\rm st}$ distance level of $T$. If $Q$ is an arbitrary maximal geodesic of $G$ that contains as a subpath the path $y-x-z$, then Lemma~\ref{lem:equal-distance} again implies, that $u$ does not belong to $Q$, a contradiction. Hence $G$ is a tree. 

We finally claim that $G$ is a subdivided star. If this is not the case, then in $T$ (which is just $G$, rooted in $u$), there exists a vertex $x$ which lies in $k^{\rm th}$ distance level of $T$, $k\ge 1$, such that $x$ has two neighbors, say $y$ and $z$, in the $(k+1)^{\rm st}$ distance level of $T$. As in the previous paragraph we now see that a maximal geodesic of $G$ that contains as a subpath the path $y-x-z$, yields a contradiction. It follows that every vertex of $T$, except maybe $u$, is of degree either $2$ or $1$. The latter is is equivalent to the fact that $T$ is a subdivided star. We hence conclude that $\gt(G) = 1$ holds if and only if $G$ is a subdivided star. 
 
Suppose now that $G$ is a an arbitrary graph that is not complete. Then there exist vertices $x, y\in V(G)$ such that $xy\notin E(G)$. But then $V(G) \setminus \{x,y\}$ form a geodesic-transversal set of $G$ and consequently, $\gt(G) \le n(G)-2$. We can hence conclude that $\gt(G) = n(G)-1$ and and only if $G$ is a complete graph of order at least $2$. 
\qed 

To conclude the section we are going to show that the geodesic-transversal problem is NP-complete. In the study of vertex-deletion problems~\cite{Yann81}, the concept of a {\em dissociation set} (see~\cite{bresar-2017, kardos-2011, tu-2020}) was considered, which was shown in~\cite{bkk11} to be the complement of a $3$-path vertex cover in any graph. Since dissociation set problem is NP-complete even when restricted to bipartite graphs~\cite{Yann81}, we infer the following. 

\begin{theorem} {\rm \cite{bkk11, Yann81}}
\label{T3PathCovNP}
	The $3$-path vertex cover problem is NP-complete for bipartite graphs.
\end{theorem}

For additional complexity results on the $3$-path vertex cover problem see~\cite{brause-2016, katrenic-2016, tsur-2019, xiao-2017}. It is clear that in bipartite graphs the $3$-path vertex cover and the $3$-geodesic transversal coincide. Thus, Theorem \ref{T3PathCovNP} can be restated as follows:

\begin{theorem}
	\label{T3GeoTraNP}
	The $3$-geodesic-transversal problem is NP-complete for bipartite graphs.
\end{theorem}

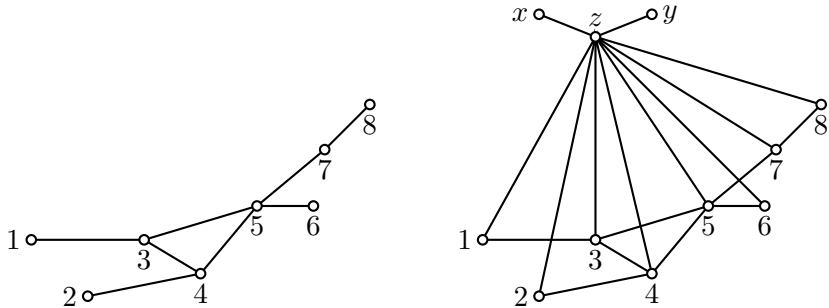
\begin{figure}[ht!]
\begin{center}
\begin{tikzpicture}[scale=1.5,style=thick,x=1cm,y=1cm]
\def\vr{1.2pt}

\begin{scope}
\coordinate(1) at (0,0);
\coordinate(2) at (0.5,-0.5);
\coordinate(3) at (1,0);
\coordinate(4) at (1.5,-0.3);
\coordinate(5) at (2,0.3);
\coordinate(6) at (2.5,0.3);
\coordinate(7) at (2.6,0.8);
\coordinate(8) at (3.0,1.2);
% \edges		
\draw (1) -- (3) -- (4) -- (2); 
\draw (4) -- (5) -- (7) -- (8); 
\draw (3) -- (5) --(6);
%  vertices
\foreach \i in {1,2,...,8}
{
\draw(\i)[fill=white] circle(\vr);
}
\draw (1) node[left] {$1$};
\draw (2) node[left] {$2$};
\draw (3) node[below] {$3$};
\draw (4) node[below] {$4$};
\draw (5) node[below] {$5$};
\draw (6) node[below] {$6$};
\draw (7) node[below] {$7$};
\draw (8) node[below] {$8$};
\end{scope}
		
\begin{scope}[xshift=4cm]
\coordinate(1) at (0,0);
\coordinate(2) at (0.5,-0.5);
\coordinate(3) at (1,0);
\coordinate(4) at (1.5,-0.3);
\coordinate(5) at (2,0.3);
\coordinate(6) at (2.5,0.3);
\coordinate(7) at (2.6,0.8);
\coordinate(8) at (3.0,1.2);
\coordinate(z) at (1.0,1.8);
\coordinate(x) at (0.5,2.0);
\coordinate(y) at (1.5,2.0);
% \edges		
\draw (1) -- (3) -- (4) -- (2); 
\draw (4) -- (5) -- (7) -- (8); 
\draw (3) -- (5) --(6);
\draw (1) -- (z) -- (2); \draw (3) -- (z) -- (4); 		
\draw (5) -- (z) -- (6); \draw (7) -- (z) -- (8); 	
\draw (x) -- (z) -- (y);	
%  vertices
\foreach \i in {1,2,...,8}
{
\draw(\i)[fill=white] circle(\vr);
}
\draw(z)[fill=white] circle(\vr);
\draw(x)[fill=white] circle(\vr);
\draw(y)[fill=white] circle(\vr);
\draw (1) node[left] {$1$};
\draw (2) node[left] {$2$};
\draw (3) node[below] {$3$};
\draw (4) node[below] {$4$};
\draw (5) node[below] {$5$};
\draw (6) node[below] {$6$};
\draw (7) node[below] {$7$};
\draw (8) node[below] {$8$};
\draw (z) node[above] {$z$};
\draw (x) node[left] {$x$};
\draw (y) node[right] {$y$};
\end{scope}

\end{tikzpicture}
\caption{A graph (left) and its reduced graph (right)}
\label{FGeoTraNP}
\end{center}
\end{figure}

Now we will prove that the geodesic-transversal problem is NP-complete for general graphs. In order to prove this, we will provide a polynomial reduction from the 3-geodesic-transversal problem to the geodesic-transversal problem. Given a graph $G$, where $V(G) = [n]=\{1,\ldots,n\}$, the reduced graph is denoted by $G'$, where $V(G') = V\cup\{x,y,z\}$ and $E(G')=E\cup\{xz,zy\}\cup\{iz:\, i \in V\}$. For an example see Fig.~\ref{FGeoTraNP}.

\begin{property} 
\label{PGeoTraNP}
	A set $S$ of vertices  is a $3$-geodesic-transversal of $G$ if and only if $S \cup \{z\}$ is a geodesic-transversal of $G'$.
\end{property}	

Property \ref{PGeoTraNP} leads to the following conclusion:
\begin{theorem}
\label{TGeoTraNP}
	The geodesic-transversal problem is NP-complete for general graphs.
\end{theorem}

%%%%%%%%%%%%%%%%%%%%%%%%%%%
%%%%%  T R E E S %%%%%%%%
%%%%%%%%%%%%%%%%%%%%%%%%%%%
%%%%%%%%%%%%%%%%%%%%%%%%%%%%%%%%%%%%
\section{The geodesic-transversal problem of trees}
\label{sec:trees}
%%%%%%%%%%%%%%%%%%%%%%%%%%%%%%%%%%%%
%Begin blue

In this section, we design an algorithm to locate a gt-set  of a tree. 

Let $T$ be a tree. A vertex of degree $1$ of a tree is a {\em leaf}. A neighbor of a leaf is a \textit{support vertex}. %A support vertex is a \textit{strong support vertex}, if it has at least two leaf neighbors, otherwise it is a \textit{weak support vertex}.  
A support vertex $u$ is an \textit{end support vertex} if $u$ is adjacent to at least $\deg(u) - 1$ leafs. 

\begin{lemma}
\label{lem:tree2}
A tree of order at least $2$ has at least one end support vertex.
\end{lemma} 

\proof
Let $T$ be a tree of order at least $2$ and let $u_1, \ldots u_k$ be the support vertices of $T$. Let $T'$ be a tree obtained from $T$ by removing all the leaves of $T$. Suppose that $\deg_{T'}(u_i) \ge 2$ for for each $i\in [k]$. Since the degree of every vertex of $T'\setminus \{u_1, \ldots, u_k\}$ is the same in $T'$ as in $T$, we would have a tree $T'$ whose every vertex is of degree at least $2$. As this is clearly not possible, there exists a vertex $u_i$ such that $\deg_{T'}(u_i) \le 1$. This in turn means that $u_i$ is an end support vertex of $T$.  
\qed

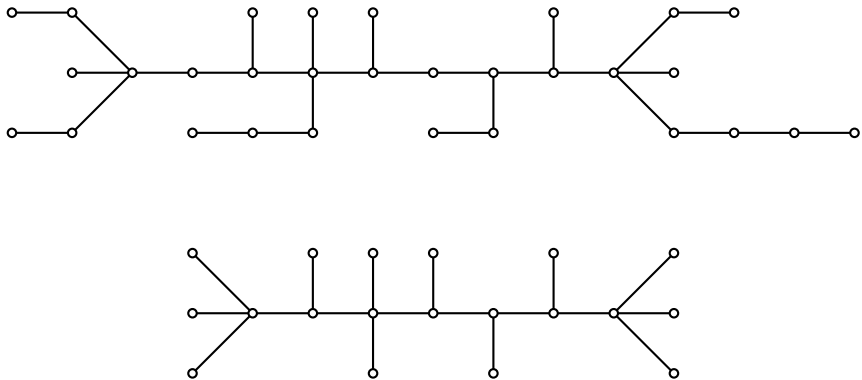
\begin{figure}[ht!]
\begin{center}
\begin{tikzpicture}[scale=0.8,style=thick,x=1cm,y=1cm]
\def\vr{2pt}

\begin{scope}
\coordinate(1) at (-1,1);
\coordinate(2) at (0,1);
\coordinate(3) at (0,0);
\coordinate(4) at (0,-1);
\coordinate(5) at (-1,-1);
\coordinate(6) at (1,0);
\coordinate(7) at (2,0);
\coordinate(8) at (3,0);
\coordinate(9) at (3,1);
\coordinate(10) at (4,1);
\coordinate(11) at (4,0);
\coordinate(12) at (4,-1);
\coordinate(13) at (3,-1);
\coordinate(14) at (2,-1);
\coordinate(15) at (5,1);
\coordinate(16) at (5,0);
\coordinate(17) at (6,0);
\coordinate(18) at (7,0);
\coordinate(19) at (7,-1);
\coordinate(20) at (6,-1);
\coordinate(21) at (8,1);
\coordinate(22) at (8,0);
\coordinate(23) at (9,0);
\coordinate(24) at (11,1);
\coordinate(25) at (10,1);
\coordinate(26) at (10,0);
\coordinate(27) at (10,-1);
\coordinate(28) at (11,-1);
\coordinate(29) at (12,-1);
\coordinate(30) at (13,-1);
% \edges		
\draw (1) -- (2) -- (6) -- (3); 
\draw (5) -- (4) -- (6) -- (7) -- (8) --(9);
\draw (8) -- (11) -- (10);
\draw (11) -- (12) -- (13) -- (14);
\draw (11) -- (16) -- (15);
\draw (16) -- (17) -- (18) -- (19) -- (20);
\draw (18) -- (22) -- (21);
\draw (22) -- (23) -- (25) -- (24);
\draw (23) -- (26);
\draw (23) -- (27) -- (28) -- (29) -- (30);
%  vertices
\foreach \i in {1,2,...,30}
{
\draw(\i)[fill=white] circle(\vr);
}
\end{scope}
		
\begin{scope}[yshift=-4cm, xshift=3cm]
\coordinate(1) at (-1,1);
\coordinate(2) at (-1,0);
\coordinate(3) at (-1,-1);
\coordinate(4) at (0,0);
\coordinate(5) at (1,0);
\coordinate(6) at (1,1);
\coordinate(7) at (2,1);
\coordinate(8) at (2,0);
\coordinate(9) at (2,-1);
\coordinate(10) at (3,1);
\coordinate(11) at (3,0);
\coordinate(12) at (4,0);
\coordinate(13) at (4,-1);
\coordinate(14) at (5,1);
\coordinate(15) at (5,0);
\coordinate(16) at (6,0);
\coordinate(17) at (7,1);
\coordinate(18) at (7,0);
\coordinate(19) at (7,-1);
% \edges		
\draw (1) -- (4) -- (2); 
\draw (3) -- (4) -- (5) -- (6);
\draw (5) -- (8) -- (7);
\draw (9) -- (8) -- (11) -- (10);
\draw (11) -- (12) -- (13);
\draw (12) -- (15) -- (14);
\draw (15) -- (16) -- (17);
\draw (18) -- (16) -- (19);
%  vertices
\foreach \i in {1,2,...,19}
{
\draw(\i)[fill=white] circle(\vr);
}
\end{scope}

\end{tikzpicture}
\caption{A tree $T$ (above) and $\SM(T)$ (below)}
	\label{FShrunkenTree}
\end{center}
\end{figure}

Let $G$ be a graph, let $v\in V(G)$ be a vertex of degree $2$, and let $x$ and $y$ be the neighbors of $u$. If $G'$ is the graph obtained from $G$ be removing the vertex $u$ and adding the edge $xy$, then we say that $G'$ is obtained from $G$ by \textit{smoothing} the  vertex $u$. Note that if the vertices $u$, $x$, and $y$ induce a triangle in $G$, then  there are two parallel edges between $x$ and $y$ in $G'$. Let further $\SM(G)$ denote a graph obtained from $G$ by smoothing all the vertices of $G$ of degree $2$. Since the smoothing operation preserves the degree of vertices, $\SM(G)$ is well-defined, that is, unique up to isomorphism. In particular, no matter in which order a smoothing of vertices of $C_n$, $n\ge 3$, is performed, we end up with $\SM(C_n) = C_2$. (The $2$-cycle $C_2$ is the graph on the vertices with two parallel edges.)  For another example see Fig.~\ref{FShrunkenTree}. 

\begin{lemma}
	\label{lem:tree1}
	If $T$ is a tree, then $\gt(T) = \gt(\SM(T))$.
\end{lemma}

\proof
Let $S$ be a gt-set of $T$. Suppose that $S$ contains a vertex $u$ with $\deg(u) = 2$. Let $P$ be the maximal path of $T$ that contains $u$ and exactly two vertices which are not of degree $2$. Such a path is indeed unique. To see it, let $x$ and $y$ be the neighbors of $u$. If ${\rm deg}(x) = 2$, then continue the path until the first vertex which is not of degree $2$ is found. Such a vertex exists since $T$ is a tree. Do the same procedure from the vertex $y$. Now, every maximal geodesic in $T$ that contains $u$, also contains $x$ and $y$. It follows that $(S\setminus \{u\})\cup \{x\}$ (or $(S\setminus \{u\})\cup \{y\}$ for that matter) is also a gt-set of $T$. Repeating this construction for every vertex of $S$ of degree $2$ we arrive at a gt-set $S'$ of $T$ which contains no vertex of degree $2$. Since $S'\subseteq V(\SM(T))$ is also a gt-set of $\SM(T)$, it follows that $\gt(\SM(T)) \le \gt(T)$. On the other hand, if $S$ is a gt-set of $\SM(T)$, then we infer that $S$ is also a gt-set of $T$, hence $\gt(T) \le \gt(\SM(T))$ also holds. 
\qed

\begin{figure}[ht!]
\begin{center}
\begin{tikzpicture}[scale=1,style=thick,x=1cm,y=1cm]
\def\vr{2pt}

\begin{scope}
\coordinate(1) at (0,0);
\coordinate(2) at (0.5,1);
\coordinate(3) at (0.9,0.1);
\coordinate(4) at (0.5,-1);
\coordinate(5) at (1.3,-2);
\coordinate(6) at (2.7,-2);
\coordinate(7) at (3.5,-1);
\coordinate(8) at (3.1,0.1);
\coordinate(9) at (4,0);
\coordinate(10) at (3.5,1);
\coordinate(11) at (2,0.8);
\coordinate(12) at (2.5,1.8);
\coordinate(13) at (1.5,1.8);
% \edges		
\draw (1) -- (3) -- (2) -- (3); 
\draw (3) -- (4) -- (5) -- (6) -- (7) -- (8) -- (11) -- (3);
\draw (13) -- (11) -- (12);
\draw (10) -- (8) -- (9);
%  vertices
\foreach \i in {1,2,...,13}
{
\draw(\i)[fill=white] circle(\vr);
}
\foreach \i in {3,5,8,11}
{
\draw(\i)[fill=black] circle(\vr);
}
\end{scope}
		
\begin{scope}[xshift=6cm, yshift=-0.7cm]
\coordinate(1) at (0.3,-0.7);
\coordinate(2) at (0.3,0.7);
\coordinate(3) at (1,0);
\coordinate(4) at (2,1);
\coordinate(5) at (1.5,2);
\coordinate(6) at (2.5,2);
\coordinate(7) at (3,0);
\coordinate(8) at (3.7,0.7);
\coordinate(9) at (3.7,-0.7);
% \edges		
\draw (1) -- (3) -- (2); 
\draw (3) -- (7) -- (4) -- (5);
\draw (5) -- (4) --(3); \draw (4) --(6); 
\draw (8) -- (7) -- (9);
%  vertices
\foreach \i in {1,2,...,9}
{
\draw(\i)[fill=white] circle(\vr);
}
\foreach \i in {3,4,7}
{
\draw(\i)[fill=black] circle(\vr);
}
\end{scope}

\end{tikzpicture}
	\caption{A graph $G$ (left) with $\gt(G) = 4$,  and $\SM(G)$ (right) with $\gt(\SM(G))=3$}
	\label{FCounterExSM}
\end{center}
\end{figure}
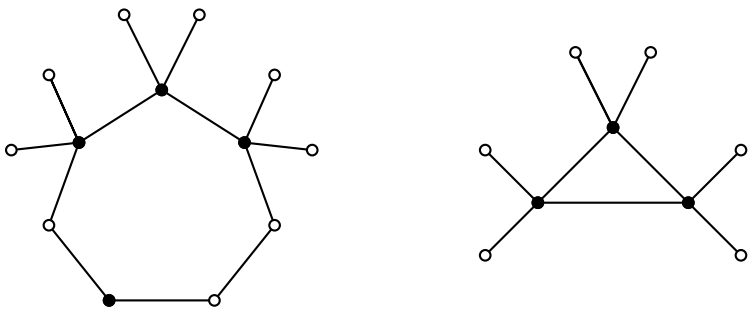

Lemma~\ref{lem:tree1} does not hold for an arbitrary graph $G$, even when $SM(G)$ does not contain parallel edges. See Fig.~\ref{FCounterExSM}, where a graph $G$ is show for which we have $\gt(G) = 4$ and $\gt(\SM(G)) = 3$.

\begin{lemma}
\label{lem:tree3}
Let $T$ be a tree with no vertices of degree $2$. Let $u$ be an end support vertex of $T$ and $u_1, \ldots, u_s$ the leaves adjacent to $u$. Then $\gt(T) = \gt(T \setminus \{u, u_1, \ldots, u_s \})+1$. Moreover, there exists a gt-set $S$ of $T$ such that $u\in S$.
\end{lemma}

\proof
Since $T$ has no vertices of degree $2$, the end support vertex $u$ is adjacent to at least two leaves, that is, $s\ge 2$. If $T$ is a star, and hence $u$ being the center of it, then the assertion of the lemma is clear. In the rest of the proof we may thus assume that $u$ has at least one non-leaf neighbor,
and since $u$ is an end support vertex, it has only one non-leaf neighbor. We denote the latter vertex by $w$, and let $T'$ be the component of $T-u$ that contains the vertex $w$. 

Let $S$ be a gt-set of $T$. Since $s\ge 2$, we see that $|S \cap \{u, u_1, \ldots, u_s \}| \ge 1$, for otherwise the geodesic $u_1, u, u_2$ would not be hit. Moreover, $|S \cap \{u, u_1, \ldots, u_s \}| = 1$. If $u_i\in S$ for some $i\in [s]$, then $(S \setminus \{u_i\})\cup \{u\}$ is also a gt-set of $T$. This proves the last assertion of the lemma and we may without loss of generality assume in the rest that $u\in S$. 

We claim now that $S\cap V(T')$ is a gt-set of $T'$. Indeed, since $\deg_{T'}(w) \ge 2$, no maximal geodesic of $T'$ can be hit by $u$. That is, only the vertices from $T'$ can be used to hit the maximal geodesics of $T'$, hence the claim. It follows that $\gt(T) = 1 + \gt(T') = 1 + \gt(T \setminus \{u, u_1, \ldots, u_s \})$
and we are done. 
\qed

Here, an algorithm is designed to construct a gt-set $S$ of an arbitrary tree $T$.

\begin{algorithm}[hbt!]
\label{al:tree}
\caption{A gt-set of a tree}
\label{alg:gt-set-tree} 
\KwIn{A tree $T$.}
\KwOut{A gt-set $S$ of $T$.}
\BlankLine
{
$S=\emptyset$;\\
$T=SM(T)$ (i.e., perform the smoothing operation on each vertex of degree $2$ in $T$).\\
\While{$|V(T)|>0$}
{
 identify an arbitrary end support vertex $p$ of $\SM(T)$;\\
 $S=S\cup\{p\}$;\\
 $T=T\setminus \{p,p_1\ldots,p_t\}$, where $p_1,\ldots,p_t$ are leaf neighbors of $p$;\\ 
 $T=\SM(T)$. \\
}
}
\end{algorithm}

\begin{theorem}
\label{thm:tree}
Given a tree $T$, Algorithm~\ref{al:tree} determines a gt-set of $T$ in linear time.  
\end{theorem}

The proof of correctness of Algorithm \ref{alg:gt-set-tree} follows from Lemmas~\ref{lem:tree2}, \ref{lem:tree1}, and \ref{lem:tree3}. The time complexity of the algorithm is clearly linear. 

To see that the smoothing operation performed in Line 2 and Line 7 of Algorithm~\ref{alg:gt-set-tree} is necessary, consider the tree $T$ in Fig.~\ref{FCounterExShTree}. Note first that $\SM(T) = 4$. Assuming that Line 2 and Line 7 would be removed from the algorithm, the modified algorithm would return a wrong value $5$. On the other hand, Algorithm~\ref{alg:gt-set-tree} first produces $\SM(T)$. Then, after two while loops (after selecting two end support vertices), another smoothing operation at Line 7 is needed. This in turn guarantees that the algorithm will end after two additional selections of end support vertices, and hence will return the correct value $4$. 

\begin{figure}[ht!]
\begin{center}
\begin{tikzpicture}[scale=0.5,style=thick,x=1cm,y=1cm]
\def\vr{4pt}

\begin{scope}
\coordinate(1) at (-1.5,1.5);
\coordinate(2) at (-1.5,0);
\coordinate(3) at (-1.5,-1.5);
\coordinate(4) at (0,0);
\coordinate(5) at (1.5,0);
\coordinate(6) at (3,0);
\coordinate(7) at (3,1.5);
\coordinate(8) at (3,3);
\coordinate(9) at (1.5,4.5);
\coordinate(10) at (3,4.5);
\coordinate(11) at (4.5,4.5);
\coordinate(12) at (3,-1.5);
\coordinate(13) at (3,-3);
\coordinate(14) at (1.5,-4.5);
\coordinate(15) at (3,-4.5);
\coordinate(16) at (4.5,-4.5);
\coordinate(17) at (4.5,0);
\coordinate(18) at (6.0,0);
\coordinate(19) at (7.5,1.5);
\coordinate(20) at (7.5,0);
\coordinate(21) at (7.5,-1.5);
% \edges		
\draw (1) -- (4) -- (2); \draw (3) -- (4);
\draw (2) -- (4) -- (5) -- (6) -- (7) --(7) -- (8) -- (9);
\draw (10) -- (8) -- (11);
\draw (6) -- (17) -- (18) -- (19);
\draw (20) -- (18) -- (21);
\draw (6) -- (12) -- (13) -- (14);
\draw (15) -- (13) -- (16);
%  vertices
\foreach \i in {1,2,...,21}
{
\draw(\i)[fill=white] circle(\vr);
}
\end{scope}

\end{tikzpicture}
	\caption{Tree $T$}
	\label{FCounterExShTree}
\end{center}
\end{figure}

%%%%%%%%%%%%%%%%%%%%%%%%%%%
%%%%%  C A C T U S %%%%%%%%
%%%%%%%%%%%%%%%%%%%%%%%%%%%
%%%%%%%%%%%%%%%%%%%%%%%%%%%%%%%%%%%%%
 \section{Fast algorithm on spread cactus graphs}
 \label{sec:cactus}
%%%%%%%%%%%%%%%%%%%%%%%%%%%%%%%%%%%%%

A connected graph in which each edge belongs to at most one cycle is a {\em a cactus graph}. We further restrict our attention to the subclass of cactus graphs in which every vertex belongs to at most one cycle, and call them {\em spread cactus graphs}. They are exactly the graphs that have neither a diamond nor a butterfly as a topological minor~\cite{nrsw-1972}. Every block in these graphs is either $K_2$ or a cycle, and cycle blocks do not intersect other cycle blocks. The blocks in a spread cactus have a tree structure, and they contain leaves or {\em leaf-cycles}, where the latter are defined as the cycle blocks, which intersect  only one $K_2$-block.

\begin{figure}[ht!]
\begin{center}
\begin{tikzpicture}[scale=0.5,style=thick]
\def\vr{4pt}
%% vertices defined %%
\path (-7,1.5) coordinate (2x);
\path (-9,0) coordinate (2y);
\path (-11,1.5) coordinate (2z);
\path (-5,3) coordinate (2a);
\path (-9,-2) coordinate (2b);
\path (-13,3) coordinate (2c);

\draw (-9,3.5) node [text=black]{$C_3(1,2,3)$}; 

%% edges %%

\draw (2x) -- (2y) -- (2z) -- (2x); 
\draw (2x) -- (2a); 
\draw (2y) -- (2b); 
\draw (2z) -- (2c); 

%% vertices %%

\draw (2x) [fill=white] circle (\vr);
\draw (2y) [fill=white] circle (\vr);
\draw (2z) [fill=white] circle (\vr);
\draw (2a) [fill=white] circle (\vr);
\draw (2b) [fill=white] circle (\vr);
\draw (2c) [fill=white] circle (\vr);

\path (2,0.5) coordinate (1x);
\path (0,-1) coordinate (1y);
\path (-2,0.5) coordinate (1z);
\path (2,2) coordinate (1a);
\path (-2,2) coordinate (1c);

%% edges %%

\draw (0,3.5) node [text=black]{$C_3(1,2)$}; 

\draw (1x) -- (1y) -- (1z) -- (1x); 
\draw (1x) -- (1a); 
\draw (1z) -- (1c); 

%% vertices %%

\draw (1x) [fill=white] circle (\vr);
\draw (1y) [fill=white] circle (\vr);
\draw (1z) [fill=white] circle (\vr);
\draw (1a) [fill=white] circle (\vr);
\draw (1c) [fill=white] circle (\vr);

\path (11,1.5) coordinate (x);
\path (9,0) coordinate (y);
\path (7,1.5) coordinate (z);
\path (9,-1.7) coordinate (b);

%% edges %%

\draw (x) -- (y) -- (z) -- (x); 
\draw (y) -- (b); 

\draw (9,3.5) node [text=black]{$C_3(1)$};
%% vertices %%

\draw (x) [fill=white] circle (\vr);
\draw (y) [fill=white] circle (\vr);
\draw (z) [fill=white] circle (\vr);
\draw (b) [fill=white] circle (\vr);

\end{tikzpicture}
\caption{Net, bull, and paw}
\label{fig:threesmall}
\end{center}
\end{figure}
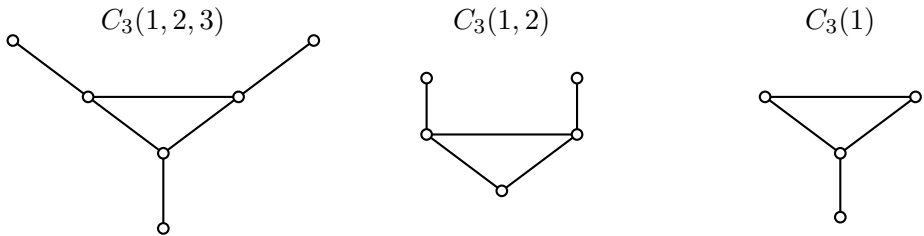

As usual, let $C_n$ denote an $n$-cycle. Let $C$ be an $n$-cycle with vertices $\{v_1,\ldots,v_n\}$, and let $I\subseteq [n]$ be a set of indices of vertices in $V(C)$. By $C_n(I)$ we denote the graph obtained from $C$ by attaching a leaf $v_i'$ to the vertex $v_i\in V(C)$ for every $i\in I$. If $I=\{i_1,\ldots,i_k\}$, then we will simplify the notation $C_n(\{i_1,\ldots,i_k\})$ to  $C_n(i_1,\ldots,i_k)$. For instance, $C_3(1,2,3)$ denotes the net graph, $C_3(1,2)$ is known as the bull graph, $C_3(1)$ is the paw graph, while $C_4(1)$ is the $P$-graph; see Fig.~\ref{fig:threesmall} for the former three graphs. 

We start our discussion by constructing an algorithm that finds a minimum geodesic transversal in the graphs $C_n(I)$ for all $n\ge 3$ and any index set $I\subseteq [n]$. Note that $C_n(I)$ are spread cactus graphs with only one cycle and no two $K_2$-blocks intersect. 

Consider $C_n(I)$, where $I=\{i_1,\ldots,i_k\}$ and $i_1<i_2<\cdots <i_k$. In the following, these indices will be taken modulo $k$. If $j\in [k]$, then we set $P^j$ to be a $v_{i_{j}}, v_{i_{j+1}}$-path along $C_n(I)$, that is, the path on vertices $v_{i_{j}}, v_{i_{j}+1},\ldots, v_{i_{j+1}}$. If $j=k$, this thus means that $P^k$ is the path on vertices $v_{i_k},v_{i_k+1},\ldots,v_1,\ldots,v_{i_1}$. 

We claim that there exists a gt-set $S$ of $C_n(I)$ such that each path $P^j$, $j\in [k]$, contains a vertex in $S$. Indeed, if 
\begin{equation*}
i_{j+1}-i_{j}\le \Big\lfloor\frac{n}{2}\Big\rfloor,
\label{eq:ij}
\end{equation*}
then $P^j$ lies on the maximal geodesic between $v_j'$ and $v_{j+1}'$. Now, if a gt-set $S$ contains $v_j'$ (resp., $v_{j+1}'$), then $S'=(S-\{v_j'\})\cup\{v_j\}$ (resp., $S'=(S-\{v_{j+1}'\})\cup\{v_{j+1}\}$) is clearly a gt-set of 
$C_n(I)$. On the other hand, if 
\begin{equation*}
i_{j+1}-i_{j}> \Big\lfloor\frac{n}{2}\Big\rfloor,
\label{eq:ij2}
\end{equation*}
then either $P^j$ contains a maximal geodesic between two vertices in $C$, or there is a maximal geodesic between $v_j'$ and $v_{j+1}$. Hence we may assume that $P^j$ contains a vertex in $S$. 

To state the next lemma, we introduce the following concept. In the graph $C_n(I)$, where $I=\{i_1,\ldots,i_k\}$, we say that $j\in [k]$ is {\em lonely}, if $i_{j+1}-i_{j-1}>\lfloor\frac{n}{2}\rfloor+1$. 

\begin{lemma}
\label{lem:gt-one-cycle}
If $n\ge 3$ and $I=\{i_1,\ldots,i_k\}$, where $0\le k\le n$, then
$$\gt(C_n(I))=\left\{ \begin{array}{ll}
\vspace*{1mm}
2; & k\le 3,\\
\vspace*{1mm}
\frac{k+1}{2}; & k\ge 5\textrm{ odd,}\\
\vspace*{1mm}
\frac{k}{2}+1; & k\ge 4\textrm{ even, and there exist lonely }j_1,j_2\in[k], j_1\textrm{ odd}, j_2\textrm{ even}, \\ 
\frac{k}{2}; & \textrm{otherwise.}
\end{array} \right.$$
\label{lem:C_n+}
\end{lemma}

\proof 
Set $G=C_n(I)$ and use the notation for vertices of $G$ as established before the lemma. Let $S$ be a gt-set of $G$. Then, as noted above, we may assume that $S\cap V(G) \subseteq C$.

We start with the case $k=|I|=0$, that is, $G=C_n$. In this case, $S=\{v_1,v_i\}$, where $i=\lfloor\frac{n}{2}\rfloor$, is clearly a gt-set of $G$, yielding $\gt(G)=2$. When $k\in\{1,2\}$, and assuming without loss of generality that $1\in I$, again the set $S=\{v_1,v_i\}$, where $i=\lfloor\frac{n}{2}\rfloor$, is a gt-set of $G$. Next, let $k=3$, and assume without loss of generality that $1\in I$. If the set $S=\{v_1,v_i\}$, where $i=\lfloor\frac{n}{2}\rfloor$, is not a gt-set of $G$, then we may assume that $1<i_2< i_3<\lfloor\frac{n}{2}\rfloor$ (the case when $\lfloor\frac{n}{2}\rfloor <i_2<i_3$ can be dealt with in a similar way).  However, then $S=\{v_2,v_2+\lfloor\frac{n}{2}\rfloor\}$ is a gt-set of $G$, yielding $\gt(G)=2$. The first line of the equality of the lemma is thus established. We next consider $k\ge 4$ and distinguish two cases.

Let $k$ be odd, $k\ge 5$. Assume that for every even $j\in [k]$, we have $i_{j+1}-i_{j-1}\le \lfloor\frac{n}{2}\rfloor+1$. Then the set $S=\{v_{i_j}:\, i_j\in I \textrm{ and } j \textrm{ odd}\}$ is a gt-set of $G$ with $|S|=\frac{k+1}{2}$. Indeed, since a maximal geodesic in $C_n$ is of length $\lfloor\frac{n}{2}\rfloor$, every maximal geodesic in $C_n(I)$ has at least one leaf as an endvertex, from which we derive that it contains a vertex $v_{i_j}$, where $j$ is odd. In the second case we may assume without loss of generality that $i_{3}-i_{1}> \lfloor\frac{n}{2}\rfloor+1$. Then $S=\{v_i:\, i=i_3-\lfloor\frac{n}{2}\rfloor-1 \textrm{ or }i>1\textrm{ odd}\}$ is a gt-set of $G$ with  $|S|=\frac{k+1}{2}$. 

Finally, let $k$ be even, $k\ge 4$. Suppose first that for every even $j\in[k]$ we have $i_{j+1}-i_{j-1}\le\lfloor\frac{n}{2}\rfloor+1$. Then, we derive in the same way as in the case of odd $k$ that the set $S=\{v_{i_j}:\, i_j\in I \textrm{ and } j \textrm{ odd}\}$ is a gt-set of $G$ with $|S|=\frac{k}{2}$. In a similar way we conclude that $\gt(G)=\frac{k}{2}$ if for every odd $j\in[k]$ we have $i_{j+1}-i_{j-1}\le\lfloor\frac{n}{2}\rfloor+1$. In the second case there exist a lonely odd $j_1\in [k]$ and a lonely even $j_2\in [k]$. Then the path $P^t$ between $v_{i_{j_t-1}}$ and $v_{i_{j_t+1}}$ is of length at least $\lfloor\frac{n}{2}\rfloor+2$, which implies that this path contains a maximal geodesic of length $\lfloor\frac{n}{2}\rfloor$, which does not involve $v_{i_{j_t+1}}$ nor $v_{i_{j_t-1}}$. Since a gt-set must hit both paths $P^{t}$, we infer that $\gt(G)>\frac{k}{2}$. It is easy to see that $\gt(G)\le \frac{k}{2}+1$ by using a similar construction as in the case when $k$ is odd. 
\qed

From the proof it is also clear that a gt-set of a graph $C_n(I)$ can be efficiently computed. If the set $I$ is a part of the input, the computation can be done in time linear in the size of $I$. 

Next, we determine a minimum geodesic transversal set $S$ in a graph $C_n(I)$ in which some of the vertices are declared in advance to be in $S$. This situation appears naturally in the construction of an algorithm for determining a gt-set of a unicyclic graph presented later. 

Let $A\subseteq [n]$ be the set of indices of the vertices of the cycle of $C_n(I)$ such that every $v_i$, $i\in A$, is predetermined to be in a geodesic transversal set $S$ of $C_n(I)$. Denote by $C_n(I,A)$ the graph $C_n(I)$ together with the requirement that vertices indexed by elements from $A$ must lie in a geodesic transversal set. The algorithm for constructing a minimum geodesic transversal of $C_n(I,A)$ is based on the constructions from the proof of Lemma~\ref{lem:C_n+}. In Algorithm~\ref{al:C_n+}, the notation of vertices $v_i\in V(C_n)$ is simplified to $i$. The indices from  $A=\{a_1,\ldots,a_t\}$ are ordered cyclically as follows: $$a_1<a_2<\cdots<a_t<a_{t+1}=a_1,$$ by which the main \texttt{while} loop is performed at least once (and is performed exactly once when $A=\{a_1\}$). The correctness of Algorithm~\ref{al:C_n+} can be proved by using similar arguments as in the proof of Lemma~\ref{lem:C_n+}.

\begin{algorithm}[hbt!]
\caption{A minimum geodesic transversal of $C_n(I,A)$}
\label{al:C_n+}
\KwIn{Cycle on $V(C_n)=\{1,\ldots,n\}$, a leaf attached to $i$, where $i\in I$, and $A\subseteq [n]$.}
\KwOut{Minimum geodesic transversal $S$ of $C_n(I)$ containing $A$.}
\BlankLine
{
	    $S=A$;\\
	    Order $A:a_1<a_2<\cdots <a_t<a_{t+1}=a_1$; \\
	    $i=1$;\\
	\While{$i\le t$}
	{
	    let $I_i=\{j\in I:\, a_i< j < a_{i+1}\}=\{j_1,\ldots,j_k\}$ and $j_0=a_i,j_{k+1}=a_{i+1}$;\\
	  \If{$k$ odd}
	  {\If{$\forall \ell\in[\frac{k+1}{2}]:\,j_{2\ell}-j_{2(\ell-1)}\le \lfloor\frac{n}{2}\rfloor+1$ }
	  {$S=S\cup\{j_{2\ell}:\,\ell\in[\frac{k-1}{2}]\}$;}
	  \Else 
	  {let $m\in [\frac{k+1}{2}]$, where $j_{2m}-j_{2m-2}>\lfloor\frac{n}{2}\rfloor+1$;\\
	   $S=S\cup\{j_{2\ell}:\,\ell\in[\frac{k-1}{2}]\}\bigcup\{j_{2m-2}+\lfloor\frac{n}{2}\rfloor+1\}$; }
	}
	\Else{let $\ell=0$;\\
	\While{$\ell\le k$}
	{ \If{$j_{\ell+2}-j_{\ell}\le \lfloor\frac{n}{2}\rfloor+1$}
	   {$S=S\bigcup \{j_{\ell+2}\}; \ell=\ell+2$;}
	  \Else
	  {$S=S\bigcup \{j_{\ell}+\lfloor\frac{n}{2}\rfloor+1\}$;\\
	  \If{$j_{\ell+1} - j_\ell \le \lfloor\frac{n}{2}\rfloor+1$}
	  {$S=S\bigcup \{j_{\ell+3}\}; \ell=\ell+3$;}	   
	  \Else 
	  {$S=S\bigcup \{j_{\ell+2}\}; \ell=\ell+2$;\\}	   
	}
	}
	}
	$i=i+1$;
}
}
\end{algorithm}

We continue by presenting an algorithm for determining a gt-set of a unicyclic graph. (This part is written mostly for intuition purposes. Algorithm~\ref{al:strongcactus} deals also with the special case when $G$ is unicyclic.)  Let $G$ be a unicyclic graph, and $C$ the cycle in $G$ of length $n$.
If $G$ is isomorphic to $C_n$, then $\gt(G)=2$. Otherwise, let $G'=G-E(C)$,  let $T_1,\ldots,T_r$ be the nontrivial components of $G'$, and let $v_1,\ldots,v_r$ be the vertices of $C$, where $v_i$ belongs to $T_i$ for all $i\in[r]$.  Clearly, each $T_i$ is a tree on at least two vertices. If $T_i$ is a path, then by the smoothing operation, and the fact that $\gt(\SM(T_i))=\gt(T_i)$, we may assume that $T_i$ is isomorphic to $P_2$, that is, $v_i$ has a leaf attached. In this case we set $S_i = \emptyset$.  Otherwise, $T_i$ has vertices of degree at least $3$, and we perform the algorithm for obtaining a gt-set $S_i$ of a tree $T_i$. It is easy to see that the sets $S_i$, $i\in[r]$, are subsets of a gt-set of $G$. There are three possibilities:
\begin{enumerate}
\item[(i)] $v_i\in S_i$;
\item[(ii)] $v_i\notin S_i$, but all neighbors of $v_i$ in $T_i$ are in $S_i$;
\item[(iii)] $v_i\notin S_i$, and there is a neighbor of $v_i$ in $T_i$ that is not  in $S_i$.
\end{enumerate} 
Turning back our attention to $G$, after {gt-sets of trees $T_i$ are obtained, the above possibilities yield different cases by which we complete the construction of a gt-set of $G$. Note that all  maximal geodesics within trees $T_i$ are hit by the sets $S_i$, hence it remains to consider the maximal geodesics that pass some vertices of $C$. The problem can be translated to determination of a minimum geodesic transversal of $C_n(I,A)$. In particular, all vertices $v_i$ that are in $S_i$ (possibility (i)) are considered to be in the set $A$, all vertices $v_i$ that are not in $S_i$ and have a neighbor in $T_i$ that is not in $S_i$ (possibility (iii)) are considered to be in $I$. Finally, the vertices $v_i\notin S_i$ for which possibility (ii) appears are in neither of the sets $A$ and $I$ (the same holds for the vertices of $C$ that are isolated in $G'$). Perform Algorithm~\ref{al:C_n+} on  $C_n(I,A)$, and let $S$ be the output of the algorithm. Finally, $S'=S\cup\bigcup_{i=1}^r{S_i}$ is a gt-set of $G$.  

We follow with two auxiliary results that will be a key for the algorithm for determining a gt-set of a spread cactus graph. We need some more notation. A vertex $v$ in a graph $G$ is {\em heavy} if $\deg_G(v)\ge 3$. Next, a heavy vertex $v$ is a {\em boundary heavy vertex} if at most one component of $G-v$ is not a path. If $v$ is a heavy vertex, then let $P^v$ denote the subset of $V(G)$ containing $v$ and every vertex of degree at most $2$ that can be reached from $v$ on a path that does not contain heavy vertices. 

\begin{lemma}
\label{lem:boundaryvertex}
If $G$ is a graph and $v$ a boundary heavy vertex in $G$ such that $G-v$ has more than two components, then $\gt(G)=1+\gt(G-P^v)$.
\end{lemma}

\proof
Since $P_t$ contains two leaves, there is a maximal geodesic that lies in $P_v$. Hence $\gt(G)\ge 1+\gt(G-P^v)$. Since every maximal geodesic in $G$ that contains a vertex in $P^v$ contains also $v$, we infer $\gt(G)=1+\gt(G-P^v)$.
\qed

Consider now a graph $G$ in which some of the vertices are declared to be in a geodesic transversal, and denote by $A_G$ the set of such vertices in $G$. (This situation appears naturally within an algorithm for determining a gt-set of $G$, where in the process of building a gt-set some of the vertices are already put in the set.) Let $C:v_1,\ldots,v_n,v_1$ be a cycle in $G$, let $A=A_G\cap V(C)$, and let $I$ be the set of vertices $v_i$, $i\in [n]$, which are adjacent to a leaf. We say that $C$ is a {\em boundary cycle} in $G$ if there exists at most one vertex $v_j\in V(C)$, where $v_j\notin I\cup A$, such that $v_i$ has a neighbor outside $C$.

\begin{lemma}
\label{lem:boundarycycle}
Let $G$ be a graph, $C$ a boundary cycle in $G$, $I$ support vertices of $C$, $A$ the set of vertices in $C$ that belong to $A_G$, and $x\in V(C)$ be adjacent to a non-leaf vertex outside $C$. Let $S_C$ be a minimum geodesic transversal of $C_n(I,A)$ and $S_C'$ a minimum geodesic transversal of $C_n(I\cup\{x\},A)$. If $|S_C'|=|S_C|$, then $S_C'$ belongs to a minimum geodesic transversal of $G$ that contains $A_G$. Otherwise, $|S_C'|=|S_C|+1$, and $S_C$ belongs to a minimum geodesic transversal of $G$ that contains $A_G$.
\end{lemma}

\proof
Clearly, $|S_C|\le |S_C'|\le|S_C|+1$. A (minimum) geodesic transversal of $G$ must hit all maximal geodesics between two vertices in $C_n(I)$.  This implies that at least $|S_C|$ vertices from $C$ need to be in a minimum geodesic transversal of $G$ that contains $A_G$. If $|S_C'|=|S_C|$, then $S_C'$ is a better choice than $S_C$, since it hits not only all the maximal geodesics that lie between two vertices in $C_n(I)$, but also all maximal geodesics that have one endvertex in $C_n(I)$. Otherwise, when $|S_C'|=|S_C|+1$, $S_C$ belongs to a minimum geodesic transversal of $G$ that contains $A_G$.
\qed

A gt-set of a path clearly consist of a single vertex, hence we may concentrate on spread cactus graphs that are not paths.  Note that for such graphs there exists a boundary heavy vertex or a boundary cycle. Hence, using Lemmas~\ref{lem:boundaryvertex} and~\ref{lem:boundarycycle}, we propose Algorithm~\ref{al:strongcactus} for determining a gt-set of a spread cactus graph.

\begin{algorithm}[hbt!]
\caption{A minimum geodesic transversal of a spread cactus graph $G$.}
\label{al:strongcactus}
\KwIn{A spread cactus graph $G$, which is not a path.}
\KwOut{Minimum geodesic transversal $S$ of $G$.}
\BlankLine
{
	    $S=\emptyset$;\\
\While{there is a heavy vertex in $G$}
	{
	\If{there is a boundary heavy vertex $v$ that lies on no cycle}
	 {$S=S\cup\{v\}$; $G=G-P^v$;}
	 \ElseIf{there is a boundary cycle $C=C_n(I,A)$, where $A=V(C)\cap S$,}
	 {
	 \If{$x$ a vertex in $C$ with a non-leaf neighbor}
	  {let $S_C$ a minimum geodesic transversal of $C_n(I,A)$ and $S_C'$ a minimum geodesic transversal of $C_n(I\cup\{x\},A)$;\\
	\If{$|S_C|=|S'_C|$} 
	   {$S=S\cup S'_C$; remove from $G$ all vertices of $C_n(I)$ and all   vertices of degree $2$ reachable by a path from $x$;  }
	\Else {$S=S\cup S_C$; $G=G-\Bigl(V(C_n(I))\setminus\{x\}\Bigr)$;}
		}
	\Else{$G=C_n(I,A)$, where $A=V(G)\cap S$, and let $S'$ be a minimum geodesic transversal of $G$ containing $A$; $S=S\cup S'$;} 
	}
	\Else{let $v$ be a boundary heavy vertex lying on a cycle;\\ 
	  \If{$\deg(v)=3$}
	     {smooth out the path $P_v$ so that $v$ is adjacent to a leaf} 
	  \Else
	    {$S=S\cup\{v\}$;\\ 
    	 $G=G-P^v$.}}
	 		
	}
}
\end{algorithm}

\begin{theorem}
\label{thm:cactus}
Given a spread cactus graph $G$, which is not a path, Algorithm~\ref{al:strongcactus} determines a gt-set of $G$ in linear time.  
\end{theorem}

\proof
By the above observations, if $G$ is a non-path spread cactus graph, then $G$ contains a heavy vertex $v$. Now, there are three possibilities: $v$ is a boundary heavy vertex that does not lie on a cycle (Line 3), $v$ lies on a cycle and its degree is at least $4$ (Line 18), or $v$ lies on a cycle and its degree is $3$. (By Line 17, $v$ can be made adjacent to a leaf.) If the latter holds for all heavy vertices of a cycle with at most one exception, then we have a boundary cycle (Line 5). The correctness of the first case and the second case (Line 3 and 18, resp.) follows from Lemma~\ref{lem:boundaryvertex}, the correctness of the second case (Lines 5-13) follows from Lemma~\ref{lem:boundarycycle}. The case when $v$ is a boundary heavy vertex with degree $3$ that lies on a cycle (Lines 16-17) follows similar  arguments as in the proof of Lemma~\ref{lem:tree1}.

An implementation of the algorithm uses a tree-like structure of a spread cactus graph, which can be obtained by a BFS search. Finding a boundary heavy vertex can be done by using a reversed order of the BFS, and all cases of the \texttt{if-then-else} condition can be checked in linear time with respect to the number of vertices that they involve. In particular, the case when there is a boundary cycle (lines 5-13) can be realized in linear time by applying Algorithm~\ref{al:C_n+} twice.
\qed

%%%%%%%%%%%%%%%%%%%%%%%%%%%%%%%%%%
%%%%%  C O N C L U D I N G %%%%%%%
%%%%%%%%%%%%%%%%%%%%%%%%%%%%%%%%%%

\section{Conclusion and future work}
\label{sec::conclusion}
A new concept of geodesic-transversal is introduced in this paper. In addition to NP-completeness, polynomial time algorithms are derived for arbitrary trees and spread cactus graphs. The potential future research is to investigate the complexity status of this problem for important interconnection networks such as butterfly networks and hypercubes, as well as for other classes of graphs such as bipartite graphs and chordal graphs. As mentioned in the initial part of the paper, it would be interesting to study how the geodesic-transversal can be used to model distance-based combinatorial problems in large-scale network analysis. 
%%%%%%%%%%%%%%%%%%%%%%%%%%%%%
\section*{Acknowledgments}
%%%%%%%%%%%%%%%%%%%%%%%%%%%%%

This work was supported and funded by Kuwait University, Research Project No.\ (QI 01/20).

%%%%%%%%%%%%%%%%%%%%%%%%%%%


\begin{thebibliography}{99}
	
\bibitem{AnSc91} 
T.~Andreae, M.~Schughart, Zs.~Tuza, 
Clique-transversal sets of line graphs and complements of line graphs, 
Discrete Math.\ 88 (1991) 11--20.

\bibitem{brkb-2017} C.~Brause, R.~Krivo\v{s}-Bellu\v{s}, 
On a relation between $k$-path partition and $k$-path vertex cover,
Discrete Appl.\ Math.\ 223 (2017) 28--38.

\bibitem{brause-2016}
C.~Brause, I.~Schiermeyer, 
Kernelization of the $3$-path vertex cover problem,
Discrete Math.\ 339 (2016) 1935--1939. 

\bibitem{bresar-2017}
B.~Bre\v{s}ar, B.~L.~Hartnell, D.~F.~Rall,
Uniformly dissociated graphs,
Ars Math.\ Contemp.\ 13 (2017) 293--306
  
\bibitem{bjk13} 
B.~Bre\v sar, M.~Jakovac, J.~Katreni\v c, G.~Semani\v sin, A.~Taranenko, 
On the vertex $k$-path cover, 
Discrete Appl.\ Math.\ 161  (2013) 1943--1949.

\bibitem{bkk11} 
B.~Bre\v sar, F.~Kardo\v s, J.~Katreni\v c, G.~Semani\v sin, 
Minimum $k$-path vertex cover, 
Discrete Appl.\ Math.\ 159 (2011) 1189--1195.

\bibitem{bkss}  
B.~Bre\v{s}ar, R. Krivo\v{s}-Bellu\v{s}, G. Semani\v{s}in, P. \v{S}parl,
On the weighted $k$-path vertex cover problem, 
Discrete Appl.\ Math.\ 177 (2014) 14--18.

\bibitem{BaNa96} 
V.~Balachandhran, P.~Nagavamsi, C.~Pandu Rangan, 
Clique transversal and clique independence on comparability graphs, 
Inform.\ Process.\ Lett.\ 58 (1996) 181--184.

\bibitem{ChCh96} 
M.-S.~Chang, Y.-H.~Chen, G.J.~Chang, J.-H.~Yan, 
Algorithmic aspects of the generalized clique-transversal problem on chordal graphs,
Discrete Appl.\ Math.\ 66 (1996) 189--203. 

\bibitem{ChSi13} 
V.~Chellappan, K.M.~Sivalingam, 
Application of entropy of centrality measures to routing in tactical wireless networks,
2013 19th IEEE Workshop on Local \& Metropolitan Area Networks (LANMAN), Brussels, 2013, 1--6. 
 
\bibitem{CoFo88} 
D.G.~Corneil, J.~Fonlupt, 
The complexity of generalized clique covering, 
Discrete Appl.\ Math.\ 22 (1988/89) 109--118. 

\bibitem{DaKr95} 
E.~Dahlhaus, J.~Kratochv\'{\i}l, P.D.~Manuel, M.~Miller, 
Parallel algorithms for generalized clique transversal problems, 
Australas.\ J.\ Combin.\ 33 (1995) 3--14. 

\bibitem{DaKr97} 
E.~Dahlhaus, J.~Kratochv\'{\i}l, P.D.~Manuel, M.~Miller, 
Transversal partitioning in balanced hypergraphs, 
Discrete Appl.\ Math.\ 79 (1997) 75--90. 

\bibitem{FoZa16} 
A.~Fornito, A.~Zalesky, E.T.~Bullmore, 
Fundamentals of Brain Network Analysis, 
Academic Press, 2016.

\bibitem{FuNu16} 
S.~Funke, A.~Nusser, S.~Storandt, 
On $k$-path covers and their applications,
Proceedings of the VLDB Endowment 7 (2014) 893--902.

\bibitem{GaJo79} 
M.R.~Garey, D.S.~Johnson, 
Computers and Intractability: A Guide to the Theory of NP-Completeness, 
Freeman, New York, 1979.

\bibitem{GoFi13} 
D.~G\'{o}mez, J.R.~Figueira, A.~Eus\'{e}bio, 
Modeling centrality measures in social network analysis using bi-criteria network flow optimization problems, 
European J.\ Oper.\ Res.\ 226 (2013) 354--365.

\bibitem{GuRa00} 
V.~Guruswami,C.P.~Rangan, 
Algorithmic aspects of clique-transversal and clique-independent sets, 
Discrete Appl.\ Math.\ 100 (2000) 183--202.

\bibitem{HoSu19} 
D.A.~Hoang, A.~Suzuki, T.~Yagita, 
Reconfiguring $k$-path vertex covers, 
Lecture Notes Comp.\ Sci.\ 12049 (2020) 133--145. 

\bibitem{HoBr95} 
M.~Hofmann, T.~Braun, G.~Carle, 
Multicast communication in large scale networks,
in: Third IEEE Workshop on the Architecture and Implementation of High Performance Communication Subsystems, Mystic, CT, USA, 1995, 147--150. 

\bibitem{JoLa10} 
K.E.~Joyce, P.J.~Laurienti, J.H.~Burdette, S.~Hayasaka, 
A new measure of centrality for brain networks, 
PLoS ONE 5 (2010) e12200. 

\bibitem{kardos-2011}
F.~Kardo\v{s}, J.~Katreni\v{c}, I.~Schiermeyer,
On computing the minimum 3-path vertex cover and dissociation number of graphs,
Theoret.\ Comput.\ Sci.\ 412 (2011) 7009--7017.

\bibitem{katrenic-2016}
J.~Katreni\v{c}, 
A faster {FPT} algorithm for $3$-path vertex cover,
Inform.\ Process.\ Lett.\ 116 (2016) 273--278. 
  
\bibitem{KoSc08} 
D.~Koschützki, F.~Schreiber, 
Centrality analysis methods for biological networks and their application to gene regulatory networks,
Gene Regul.\ Syst.\ Bio.\ 2 (2008) 193--201. 

\bibitem{KwCh19} 
H.~Kwon, Y-H.~Choi, J.-M.~Lee, 
A physarum centrality measure of the human brain network, 
Sci.\ Reports 9 (2019) Article 5907. 

\bibitem{LiLi15} 
C.~Li, Q.~Li, P.~Van Mieghem, H.E.~Stanley, H.~Wang, 
Correlation between centrality metrics and their application to the opinion model, 
European Phys.\ J.\ B 88 (2015) Article 65.

\bibitem{LiLi19} 
W.~Liu, X.~Li, T.~Liu, B.~Liu, 
Approximating betweenness centrality to identify key nodes in a weighted urban complex transportation network, 
J.\ Adv.\ Transp.\ (2019) Article 9024745. 

\bibitem{Manu18} 
P.~Manuel, 
Revisiting path-type covering and partitioning problems, 
arXiv:1807.10613 [math.CO] (25 July 2018).

\bibitem{Manu19} P.~Manuel, 
On the isometric path partition problem, Discussiones Mathematicae Graph Theory, pp 1-13, 2019, DOI: https://doi.org/10.7151/dmgt.2236.

\bibitem{MoMa20} 
D.B.~Mokeev, D.S.~Malyshev, 
A polynomial-time algorithm of finding a minimum $k$-path vertex cover and a maximum $k$-path packing in some graphs,
Opt.\ Lett.\ 14 (2020) 1317--1322. 

\bibitem{nrsw-1972} 
E.A.~Nordhaus, R.D.~Ringeisen, B.M.~Stewart, A.T.~White, 
A Kuratowski-type theorem for the maximum genus of a graph, 
J.\ Combin.\ Theory Ser.\ B 12 (1972) 260--267.

\bibitem{PeGe16} 
C.~Perez, R.~Germon, 
Chapter 7 - Graph Creation and Analysis for Linking Actors: Application to Social Data,
in: R.~Layton, P.A.~Watters (eds.), Automating Open Source Intelligence, Syngress, 2016, 103--129. 

\bibitem{Prak12} 
B.A.~Prakash, 
Propagation and immunization in large networks, 
ACM Mag.\ Stud.\ 19 (2012) 56--59.

\bibitem{Rack02} 
H.~Racke, 
Minimizing congestion in general networks, 
in: The 43rd Annual IEEE Symposium on Foundations of Computer Science, 2002, Vancouver, BC, 2002, 43--52. 

\bibitem{Ren15} 
Y.~Ren, 
Betweenness Centrality and Its Applications from Modeling Traffic Flows to Network Community Detection,
PhD Thesis, University of Notre Dame, USA, 2015.

\bibitem{Sing05} 
A.~Singh, 
Load-balanced routing in interconnection networks, 
PhD thesis, Stanford University, USA, 2005.

\bibitem{StKa17} 
M.~Stewart, R.~Kannan, A.~Dvir, B.~Krishnamachari, 
CASPaR: Congestion avoidance shortest path routing for delay tolerant networks, 
Int.\ J.\ Dist.\ Sensor Networks 13(11) (2017) 1--15.

\bibitem{tsur-2019} 
D.~Tsur, 
Parameterized algorithm for $3$-path vertex cover,
Theoret.\ Comput.\ Sci.\ 783 (2019) 1--8. 

\bibitem{TuYa13} 
J.~Tu, F.~Yang, 
The vertex cover $P_3$ problem in cubic graphs, 
Inform.\ Process.\ Lett.\ 113 (2013) 481--485.

\bibitem{tu-2020}
J.~Tu, Z.~Zhang Y.~Shi,
The maximum number of maximum dissociation sets in trees, 
J.\ Graph Theory (2020) doi.org/10.1002/jgt.22627. 

\bibitem{TuZh11} 
J~Tu, W.~Zhou, A primal–dual approximation algorithm for the vertex cover $P_3$ problem, Theoret.\ Comput.\ Sci.\ 412 (2011) 7044--7048. 

\bibitem{WaZh16} 
X.~Wang, W.~Ni, K.~Zheng, R.P.~Liu, X.~Niu, 
Virus propagation modeling and convergence analysis in large-scale networks,
in: IEEE Trans.\ Inform.\ Forensics Sec.\ 11 (2016) 2241--2254.

\bibitem{XiHa16} 
W.~Xi, J.~Han, K.~Li, Z.~Jiang, H.~Ding, 
Chapter 13 - Location Inferring in Internet of Things and Big Data, 
in: R.~Buyya, R.N.~Calheiros, A.V.~Dastjerdi (eds.), Big Data, Morgan Kaufmann, 2016,  309--335.

\bibitem{xiao-2017}
M.~Xiao, S.~Kou, 
Exact algorithms for the maximum dissociation set and minimum $3$-path vertex cover problems,
Theoret.\ Comput.\ Sci.\ 657 (2017) 86--97. 

\bibitem{Xu02} 
J.-M.~Xu, 
Topological Structure and Analysis of Interconnection Networks, 
Springer, 2002.

\bibitem{XuHu19} 
Z.~Xu, X.~Huang, Y.~Deng, 
Load-balanced routing for nested interconnection networks, 
arXiv:1909.06497v2 [cs.NI] (25 Dec 2019). 

\bibitem{Yann81} 
M.~Yannakakis, 
Node-deletion problems on bipartite graphs, 
SIAM J.\ Comput.\ 10 (1981) 310--327.

\bibitem{YuGu15} 
S.~Yu, G.~Gu, A.~Barnawi, S.~Guo, I.~Stojmenovic, 
Malware propagation in large-scale networks,
IEEE Trans.\ Know.\ Data Eng.\ 27 (2015) 170--179.

\bibitem{ZhNa15} 
P.~Zhao, S.M.~Nackman, C.K.~Law, 
On the application of betweenness centrality in chemical network analysis: Computational diagnostics and model reduction, 
Combust.\ Flame 162 (2015) 2991--2998.

\bibitem{zak} 
A. \.{Z}ak, 
Generalized transversals, generalized vertex covers and node-fault-tolerance in graphs, Discrete Appl.\ Math.\ 255 (2019) 299--306.
	
\end{thebibliography}
\end{document}